# On Primes In Short Intervals
N. A. Carella, December, 2008.

**Abstract:** This note discusses the existence of prime numbers in short intervals. An unconditional elementary argument seems to prove the existence of primes in the short interval $[x, x + y]$ with $y \geq x^{1/2}(\log x)^\varepsilon$, $\varepsilon > 0$, and a sufficiently large number $x > 0$. Further, an extension of Bertrand's postulate to arithmetic progressions will be considered.

## 1 INTRODUCTION
This note discusses the existence of prime numbers in short intervals. An unconditional elementary argument seems to prove the existence of primes in the short interval $[x, x + y]$ with $y \geq x^{1/2}(\log x)^\varepsilon$, $\varepsilon > 0$, and a sufficiently large number $x > 0$. Further, an extension of Bertrand's postulate to arithmetic progression will be considered.

Currently the sharpest estimates for primes on short intervals are the followings.
(i) Unconditionally, there exists a prime in the short interval $[x, x + y]$ with $y \geq x^{.525}$, see Theorem 52.
(ii) Assuming the Riemann hypothesis, there exists a prime in the short interval $[x, x + y]$ with $y \geq x^{1/2}(\log x)^2$, see Theorem 53.

The new results contributed to the literature are the followings.

**Theorem 1.** For all sufficiently large numbers $x > 0$, and $y \geq x^{1/2}(\log x)^\varepsilon$, $\varepsilon > 0$, the interval $[x, x + y]$ contains some primes.

**Theorem 2.** Let $q \leq (\log x)^A$, $A > 0$, and let $x > 1$ be a sufficiently large number. Then the interval $[x, 2x]$ contains a prime number from the arithmetic progression $p \equiv a \bmod q$, $\gcd(a, q) = 1$.

The proofs of these results appear on page 22 and page 25 respectively. As a consequence, the first result leads to an improvement on the conditional prime gap from $p_{n+1} - p_n \leq c p_n^{1/2}(\log p_n)^2$ to the unconditional prime gap $p_{n+1} - p_n \leq c p_n^{1/2}(\log p_n)^\varepsilon$, $c > 0$ constant, see Theorems 53.

The articles [HB], [IC] and [PT] present surveys on the state of knowledge in the theory of prime numbers and related topics. Other references for other specialized subjects are given in throughout the paper.

## 2 SURVEY OF THE ZETA AND *L*-FUNCTIONS
A short survey of the properties of the zeta function and the Dirichlet *L*-functions is considered here. It will also help to establish some of the notations used throughout the paper.

**Zeta Function**
The zeta function is defined by the power series and product

$$\zeta(s) = \sum_{n=1}^{\infty} \frac{1}{n^s} = \prod_{p \geq 2} (1 - 1/p^s)^{-1}, \qquad (1)$$

where $s \in \mathbb{C}$, see [EU, p. 228], and the inverse is defined by

$$\frac{1}{\zeta(s)} = \sum_{n=1}^{\infty} \frac{\mu(n)}{n^s} = \prod_{p \geq 2} (1 + 1/p^s)^{-1}. \qquad (2)$$

Both of these power series are absolutely convergent on the complex half plane $\Re(s) > 1$. Each term $\zeta_p(s) = (1 - 1/p^s)^{-1} = \sum_{n \geq 1} p^{-ns}$ in the product is a local zeta function at *p*.

The Weierstrass product representation of the zeta function explicitly shows the known locations of its zeros. This is given by

$$\zeta(s) = \frac{1}{2}\left(\frac{2\pi}{e}\right)^s \prod_{n=1}^{\infty}\left(1 + \frac{s}{2n}\right)e^{-s/2n} \prod_{\rho}\left(1 - \frac{s}{\rho}\right)e^{s/\rho}, \qquad (3)$$

where $s \in \mathbb{C}$, and $\rho$ ranges over the zeros of $\zeta(s)$ on the critical strip $0 \leq \Re(s) \leq 1$. A simple scaling of the original definition of $\zeta(s)$ realizes a simple and useful analytic continuation on the half plane $\Re(s) > 0$.

***Theorem 3.*** (i) The analytic continuation $\zeta(s) = (1 - 2^{1-s})^{-1} \sum_{n \geq 1} (-1)^{n-1} n^{-s}$ of $\zeta(s)$ is convergent for $\Re(s) > 0$, absolutely convergent for $\Re(s) > 1$, and has a pole at $s = 1$.
(ii) The function $\zeta(s) = s/(s-1) - s\int_{1}^{\infty} \{t\} t^{-s-1} dt$ is an analytic continuation of $\zeta(s)$ for $\Re(s) > 0$ and has a pole at $s = 1$, where $\{x\} = x - [x]$ is the fractional part of *x*.

The derivation of (i) is simply $\sum_{n \geq 1}(-1)^{n-1} n^{-s} = \sum_{n \geq 1} n^{-s} - 2\sum_{n \geq 1}(2n)^{-s} = \zeta(s) - 2^{1-s}\zeta(s)$. To verify the convergent statement for $\Re(s) > 0$, let $S_{2n} = \sum_{1 \leq d \leq n}(-1)^{d-1} d^{-s}$ and $s = \sigma + it$. Now use the fact that $(2d-1)^{-s} - (2d)^{-s} = s\int_{2d-1}^{2d} x^{-s-1} dx$ to obtain $|S_{2n}| \leq |s| \sum_{1 \leq d \leq n} (2d-1)^{-\sigma-1}$. These information imply that $\left|\sum_{n \geq 1}(-1)^{n-1} n^{-s}\right| \leq |s| \sum_{d \geq 1}(2d-1)^{-\sigma-1}$, so the power series is convergent but not absolutely convergent. A more general technique of analytical continuation will be introduced later.

The closely related alternating zeta function $\zeta^*(s) = (1 - 2^{1-s})\zeta(s)$ has many properties in common with its counterpart. It is holomorphic on the complex plane, it has no poles, but it has trivial zeros at $s = -2n$, $n \geq 1$, and nontrivial zeros $s = 1 + i2\pi n/\log 2$, $n \in \mathbb{Z}$, on the line $\Re(s) = 1$. These zeros arise from the factor $1 - 2^{1-s}$.

The there are other complicated analytic continuations of $\zeta(s)$ to the punctured complex plane $\mathbb{C} - \{\,1\,\}$. Some of these are

(i) The Euler-Maclaurin expansion

$$\zeta(s) = 1 + \frac{1}{s-1} + \sum_{k=2}^{n}(-1)^k \binom{-s}{k}B_k - \frac{(-1)^n}{n}\binom{-s}{n-1}\int_1^\infty B_n(\{x\})\frac{dx}{x^{s+n}} \text{ for } \Re e(s) > 1-n. \quad (4)$$

(ii) The Hasse series

$$\zeta(s) = (1-2^{1-s})^{-1}\sum_{n=0}^{\infty}\frac{1}{2^{n+1}}\sum_{k=0}^{n}(-1)^k\binom{n}{k}\frac{1}{(k+1)^s}, \; 1 \neq s \in \mathbb{C}. \quad (5)$$

Some numerical calculations (and computations of some special values) can be implemented with these rapidly convergent series.

**Functional Equation**

The functional equation $F(x) = F(x - t)$ of a function $f(x)$ is concerned with the fixed points $t$ of the translation action $x \to x - t$ for all $x \in \mathbb{C}$. An important case is the function

$$\xi(s) = (s/2)(s-1)\pi^{-s/2}\Gamma(s/2)\zeta(s), \quad (6)$$

which has the Weierstrass product representation

$$\xi(s) = e^{A+Bs}\prod_{\rho}(1-s/\rho)e^{s/\rho}, \quad (7)$$

where $s \in \mathbb{C}$, $\rho$ ranges over the zeros of $\zeta(s)$ and $A$ and $B$ are constants.

***Theorem* 4.** (Riemann 1859) The function $\xi(s)$ is entire on the complex plane $\mathbb{C}$ and satisfies the functional equation $\xi(s) = \xi(1-s)$ for any complex number $s \in \mathbb{C}$, $\text{Im}(s) \neq 0$.

The statement $\text{Im}(s) \neq 0$ stems from the fact that $\zeta(2n) \neq \zeta(1-2n)$ for all integers $n \geq 1$. Three other important observations on the function $\xi(s)$ are the followings.
(i) The critical line $\Re e(\rho) = 1/2$ is the line of symmetry of the functional equation.
(ii) The subset of trivial zeros $\{\,-2n : n \in \mathbb{N}\,\}$ and the subset of nontrivial zeros $\{\,\rho = \beta + it : 0 \leq \beta \leq 1\,\}$ of the entire function $\xi(s)$ coincides with the zeros of the zeta function $\zeta(s)$.
(iii) The poles of the gamma function at the negative even integers cancel the trivial zeros of the zeta function at the negative even integers. Thus the trivial zeros of the zeta function are simple zeros.

Quite often the simpler functional equation

$$\pi^{-s/2}\Gamma(s/2)\zeta(s) = \pi^{(s-1)/2}\Gamma((1-s)/2)\zeta(1-s) \quad (8)$$

is used instead of $\xi(s) = \xi(1-s)$.

**Conjecture 5.** (Riemann 1859) The nontrivial zeros of the zeta function are on the critical line $\Re e(\rho) = 1/2$.

As of 2008, every known zeros of $\zeta(s)$ lies on the critical line, see [GN] for extensive computational details. The imaginary part $t \geq 0$ of the first 10 zeros $\rho_n = 1/2 + it_n$ at 12 decimal places accuracy are listed in the table.

| n | $t_n$ | n | $t_n$ |
|---|---|---|---|
| 1 | 14.1347251417346 | 6 | 37.5861781588256 |
| 2 | 21.0220396387715 | 7 | 40.9187190121474 |
| 3 | 25.0108575801456 | 8 | 43.3270732809149 |
| 4 | 30.4248761258595 | 9 | 48.0051508811671 |
| 5 | 32.9350615877391 | 10 | 49.7738324776723 |

The behavior of the zeta function on the critical strip $0 \leq \Re(s) \leq 1$ is predicted to vary as a function of the imaginary part $\text{Im}(s) = t$ of the complex number $s = \sigma + it$. The first result below is proven via the convexity principle or the Phragmen-Lindelof principle, and the second is an open problem.

**Theorem 6.** For $\varepsilon > 0$, $\zeta(1/2 + it) = O(t^{1/4+\varepsilon})$.

**Lindelof hypothesis 7.** For $\varepsilon > 0$, $\zeta(1/2 + it) = O(t^\varepsilon)$.

**Modular Theta Functions**
The modular theta functions are used to derive the analytic continuations and the functional equations of the zeta function and *L*-functions. A few other techniques are also known for accomplishing these tasks. For example, complex analysis (residue integrals). However, the proofs are longer than the ones based on theta functions.

The modular theta functions are defined by

(i) $\theta(t) = \sum_{n=-\infty}^{\infty} e^{-\pi n^2 t}$,  (ii) $\theta_0(t, \chi) = \sum_{n=-\infty}^{\infty} \chi(n) e^{-\pi n^2 t}$,  (iii) $\theta_1(t, \chi) = \sum_{n=-\infty}^{\infty} n\chi(n) e^{-\pi n^2 t}$,

where $\chi$ is a character modulo $q$ (the last two formulas are for even and odd character respectively). These are modular functions of weight 1/2 on the modular subgroup $\Gamma(2) = \{ z \to -z^{-1}, z \to z+2 \}$.

**Theorem 8.** Let $t > 0$. Then the inversion formulas are
(i) $\theta(t + i2) = \theta(t)$,                                 (ii) $\theta(t) = t^{-1/2} \theta(1/t)$,
(iii) $\theta_0(t, \chi) = q^{-1/2} \tau(\chi) t^{-1/2} \theta_0(1/t, \overline{\chi})$,   (iv) $\theta_1(t, \chi) = q^{-1/2} \tau(\chi) t^{-3/2} \theta_1(1/t, \overline{\chi})$,

for even and odd character $\chi$ respectively.

The proofs of the functional equation of the zeta function and analytic continuation are considered here. It is based on the evaluations of the Mellin transform of the theta function. For $\Re(s) > 1$, the Mellin transform of the theta function is evaluated in two ways. The first evaluation is

$$2^{-1} \int_0^\infty (\theta(t) - 1) t^{s/2-1} dt = \int_0^\infty \left( \sum_{n=1}^\infty e^{-\pi n^2 t} \right) t^{s/2-1} dt = \frac{1}{\pi^{s/2}} \sum_{n=1}^\infty \frac{1}{n^s} \int_0^\infty x^{s/2-1} e^{-x} dx = \pi^{-s/2} \Gamma(s/2) \zeta(s). \tag{9}$$

Here the change of variable $x = \pi n^2 t$ and the absolute convergence of both the integral, and the sum are used. The second evaluation is

$$2^{-1} \int_0^\infty (\theta(t) - 1) t^{s/2-1} dt = 2^{-1} \int_0^1 (\theta(t) - 1) t^{s/2-1} dt + 2^{-1} \int_1^\infty (\theta(t) - 1) t^{s/2-1} dt, \tag{10}$$

where the change of variable $x = 1/t$ leads to

$$2^{-1}\int_0^1 (\theta(t)-1)t^{s/2-1}dt = -\frac{1}{s} + 2^{-1}\int_0^1 \theta(t)t^{s/2-1}dt = -\frac{1}{s} + 2^{-1}\int_1^\infty \theta(x^{-1})x^{-s/2-1}dx$$

$$= -\frac{1}{s} + \frac{1}{1-s} + 2^{-1}\int_1^\infty (\theta(t)-1)x^{(1-s)/2-1}dx.$$
(11)

Replacing (11) back into (10) yields

$$2^{-1}\int_0^\infty (\theta(t)-1)t^{s/2-1}dt = -\frac{1}{s} + \frac{1}{s-1} + 2^{-1}\int_1^\infty (\theta(t)-1)(x^{s/2-1} + x^{(1-s)/2-1})dx.$$
(12)

Since the integral converges for all complex number $s \in \mathbb{C}$, the analytic continuation of $\zeta(s)$ holds for $s$. Moreover, the last expression is invariant under the action of the involution $s \to 1-s$, so it follows that

$$\pi^{-(s-1)/2}\Gamma((1-s)/2)\zeta(1-s) = \pi^{-s/2}\Gamma(s/2)\zeta(s).$$
(13)

**Explicit Formulae**
The explicit formula was introduced by Riemann in 1859 and proved by von Mangoldt in 1895. Today explicit formulas are standard tool in prime number theorems and related topics.

***Theorem 9.*** Let $\rho = \beta + it$ denotes the zeroes of the zeta function. Then

$$\psi(x) = x - \sum_\rho \frac{x^\rho}{\rho} - \log 2\pi - \frac{1}{2}\log(1 - 1/x^2),$$
(14)

where $x \neq p^n$.

Nowadays it is quite easy to derive this formula using the Perron integral and logarithm derivative or some other methods. The Perron integral and the logarithm derivative of the zeta function are

$$\frac{1}{i2\pi}\int_{c-i\infty}^{c+i\infty} x^s \frac{ds}{s} = \begin{cases} 0, & 0 < x < 1, \\ 1/2, & x = 1, \\ 1, & x > 1, \end{cases} \quad \text{and} \quad -\frac{\zeta'(s)}{\zeta(s)} = \sum_{n=1}^\infty \frac{\Lambda(n)}{n^s},$$
(15)

respectively, see Section 7 for a proof. The arithmetic function $\Lambda(n) = \log p$ if $n = p^m$, $p$ prime and $m \geq 1$, otherwise, $\Lambda(n) = 0$. The meromorphic function $f(s) = -\frac{\zeta'(s)}{\zeta(s)}\frac{x^s}{s}$ has all its poles at $s = 0, 1$ and at the zeros $\rho = \sigma + it \in \mathbb{C}$ of the zeta function. The goal is to evaluate integral transform of the function $f(s)$ in two ways: In terms of residues, this is

$$\sum_{\text{residues of } f} res(f) = \frac{1}{i2\pi}\int_{a-i\infty}^{a+i\infty} -\frac{\zeta'(s)}{\zeta(s)}\frac{x^s}{s}ds = x - \sum_{\rho \neq 0,1} \frac{x^\rho}{\rho} - \frac{\zeta'(0)}{\zeta(0)},$$
(16)

where the index ρ runs over the zeros. To complete the integral evaluation, use $\zeta(0) = -1/2$, $\zeta'(0) = -(\log 2\pi)/2$, and $\sum_\rho x^\rho / \rho = \log(1 - x^{-2})/2$, where the zeros $\rho = -2n$, $n \geq 1$, are trivial. On the other hand, in terms of an arithmetic function, the evaluation is

$$\frac{1}{i2\pi} \int_{a-i\infty}^{a+i\infty} -\frac{\zeta'(s)}{\zeta(s)} \frac{x^s}{s} ds = \frac{1}{i2\pi} \int_{a-i\infty}^{a+i\infty} \sum_{n=1}^{\infty} \frac{\Lambda(n)}{n^s} \frac{x^s}{s} ds = \sum_{n=1}^{\infty} \Lambda(n) \left( \frac{1}{i2\pi} \int_{a-i\infty}^{a+i\infty} \left(\frac{x}{n}\right)^s \frac{ds}{s} \right) = \sum_{n \leq x} \Lambda(n) = \psi(x) \quad (17)$$

for $x \in \mathbb{R} - \{$ prime power $p^n \}$. For other details, see [EL, p. 166], [ES, p. 58], [MO, p. 147], etc., and other information on explicit formulas for zeros of the zeta function, see [EL, p 162] and similar sources.

The nonvanishing of the zeta function on the line $\Re(s) = 1$ was established by both DelaValle Poussin and Hadamard in 1896, but the most common and easiest proof appears to be due to Mertens.

***Theorem* 10.** The zeta function $\zeta(1 + it) \neq 0$ for all $t \in \mathbb{R}$.

Proof (Mertens 1898): $\mathrm{Re}(\log \zeta(s)) = \sum_{p \geq 2} \sum_{n=1}^{\infty} \frac{1}{np^{n\sigma}} \cos(t \log p^n)$. From the cosine identity

$$3 + 4\cos\theta + \cos 2\theta = 2(1 + \cos\theta)^2 \geq 0,$$

it follows that

$$3\log\zeta(\sigma) + 4\mathrm{Re}\log\zeta(\sigma + it) + \mathrm{Re}\log\zeta(\sigma + i2t) \geq 0.$$

Reversing the logarithm gives $\zeta(\sigma)^3 \zeta(\sigma + it)^4 \zeta(\sigma + i2t) \geq 1$. Next taking the limit $\sigma \to 1$, and using the assumption $\zeta(1 + it) = 0$ gives a contradiction. ∎

Another elementary proof of the nonvanishing of the zeta function can be realized using the class number formula for quadratic fields. A quantitative version of the previous nonvanishing result is given below, it is used to effectively derive an improve version of the prime number theorem.

***Theorem* 11.** (DelaValle Poussin 1898) $\zeta(s) \neq 0$ for all $s = \sigma + it$ such that $\Re(s) \geq 1 - c/\log(|t| + 2)$, $c > 0$ constant.

**Zeros Statistics**
A few elementary results on the zeros of the zeta function are provided in this subsection. A basic fact is the location of the zeros on the complex plane $\mathbb{C}$. In the critical strip $0 < \Re(s) < 1$ all the zero are complex. This fact is extracted from the values of the power series $\sum_{n \geq 1} (-1)^{n-1} n^{-s} = (1 - 2^{1-s})\zeta(s) = 1 - 2^{-s} + 3^{-s} - \cdots > 0$ for real $s$ in the range $0 < s < 1$ and the value $\zeta(0) = -1/2$.

A well known result on the zeros of the zeta function is their density on a finite rectangle

$$R(T) = \{ s \in \mathbb{C} : -T \leq \mathrm{Im}(s) \leq T, -1 \leq \Re(s) \leq 1 \}.$$

***Theorem* 12.** (i) Let $N(T) = \#\{ s \in R(T) : \zeta(s) = 0 \}$ be the number of complex zeros on the critical strip. Then

$$N(T) = \frac{T}{2\pi} \log \frac{T}{2\pi} - \frac{T}{2\pi} + O(\log T) + S(T),$$

where the term $S(T) = \pi^{-1} \arg \zeta(1/2 + it) = O(\log T)$ for all $T \in \mathbb{R}$.
(ii) $N(T+1) - N(T) = \sum_{T < \gamma \leq T+1} 1 = O(\log T)$ for $\rho = \beta + i\gamma$.
For a proof (about two pages) based on the Cauchy integral (winding number)

$$\frac{1}{i2\pi} \int_R \frac{\xi'(s)}{\xi(s)} ds, \tag{18}$$

over the rectangle $-1 - iT$, $-1 + iT$, $2 + iT$, $2 - iT$, see [KA, p. 26].

A result for the number of zeros in a short interval states the following.

***Theorem* 13.**  (Selberg 1942)  Let $S \geq T^{1/2 + \varepsilon}$, $\varepsilon > 0$. Then

$$N(T+S) - N(T) \geq (c/2\pi) S \log T,$$

where $c > 0$ is a constant.

This result was improved by Moser in 1976 and other authors (confer the literature) to
(i) $N(T+S) - N(T) > 0$,  for $S \geq T^{1/6} \log^2 T$,
(ii) $N(T+S) - N(T) \geq cS$,  for $S \geq T^{5/12} \log^3 T$.

A *density theorem* provides an upper estimate for the cardinality of the set $N(\sigma, T) = \#\{ \rho = \sigma + it : \sigma > 1/2 \}$. This is something of the form $N(\sigma, T) \leq cT^{f(\sigma)(1-\sigma)} \log^B T$, where $B, c > 0$ are constants and $2 \leq f(\sigma)$ is a function, $1/2 \leq \sigma \leq 1$. The measure $N(\sigma, T)$ tracks the density of the zeros of $\zeta(s)$ on the strip $1/2 \leq \Re(s) \leq 1$. It is believed that $N(\sigma, T) = 0$.

***Theorem* 14.**  (Selberg 1946)  Let $\sigma \geq 1/2$ and let $T > 0$. Then

$$N(\sigma, T) = cT^{1-(\sigma-1/2)/4} \log T,$$

where $c > 0$ is a constant.

***Theorem* 15.**  ([HL])  For $\sigma \geq 1/2$ and let $T > 0$. Then

$$N(\sigma, T) < T^{12(1-\sigma)/5} \log^B T,$$

where $B > 0$ is a constant.

**Zeros Detection**. For all real number $t \in \mathbb{R}$, The real-valued function $Z(t) = e^{i\theta(t)} \zeta(1/2 + it)$ has the same zeros lying on the critical line $\Re(s) = 1/2$ as the zeta function. The angle is defined by

$$\theta(t) = -i \ln\left(\pi^{-it/2} \Gamma(1/4 + it/2) / |\Gamma(1/4 + it/2)|\right).$$

The simpler asymptotic expression $\theta(t) = \frac{t}{2}\log\frac{t}{2\pi} - \frac{t}{2} - \frac{\pi}{8} + O(1/t)$ can be derived using Sterling's formula.

A basic rule of zeros detection is the Gram's law. This law claims that there is a zero in the interval $[g_n, g_{n+1}]$, where the point $g_n$ satisfies $\theta(g_n) = n\pi$. This is due to the relation $\zeta(1/2 + it) = Z(t)\cos\theta(t) + iZ(t)\sin\theta(t) = (-1)^n Z(t)$ at the point $t = g_n$. This relation specifies at least one change of sign in the real values of the zeta function. But this law is known to fail infinitely often, consult the literature for finer details.

## Special Values

The determination of the set of special values $\{ \zeta(2n + 1) : n \in \mathbb{N} \}$ in closed forms and the set of critical zeros $\{ \beta + it : \zeta(\beta + it) = 0 \text{ and } 0 \leq \beta \leq 1 \}$ are some of the most important open problems in the theory of the zeta function.

***Theorem* 16.** (Euler) For $n \geq 1$, the values of the zeta function at the integer arguments are the followings:
(i) $\xi(-2n) = 0$,
(ii) $\zeta(1 - 2n) = -B_{2n}/(2n)$,
(iii) $\zeta(2n) = (-1)^{n-1} 2^{2n-1} \pi^{2n} B_{2n}/(2n)!$,
(iv) $\xi(2n + 1) = \text{Unknown}$,
where $B_n$ is the $n$th Bernoulli number.

Proof: For (i) use the fact that the gamma function $\Gamma(s)$ has a simple zero at $s = -n$, $n \geq 1$, and the formula $\zeta(1-s) = 2(2\pi)^{-s}\Gamma(s/2)\cos(\pi s/2)\zeta(s)$ to determine the simple zeros of $\zeta(s)$. For (ii) and (iii) evaluate the Fourier series $B_n(x) = -n!\sum_{d \geq 1}(i2\pi d)^{-n} e^{-i2\pi dx}$, $0 \leq x < 1$, see Section 7. ∎

***Conjecture* 17.** (Siegel 1966) The rational numbers $2n\zeta(1-2n) = -B_{2n}$ are uniformly distributed modulo $p$, $\gcd(n, p) = 1$.

A similar claim can be made for some special values of the L-functions attached to a character $\chi$ modulo $q$: The rational numbers $nL(1-n, \chi) = -B_{n,\chi}$ are uniformly distributed modulo $p$, $\gcd(n, p) = 1$, $\gcd(n, q) = 1$.

## Multiplicative Characters

A multiplicative character $\chi$ is a periodic, complex valued and completely multiplicative function $\chi : \mathbb{Z} \to \mathbb{C}$ on the integers. For each $q \in \mathbb{N}$, the set of characters $\hat{G} = \{1 = \chi_0, \chi_1, ..., \chi_{\varphi(q)-1}\}$ is a group of order $\varphi(q) = \prod_{p|q}(1 - 1/p)$, where $p$ ranges over the prime divisors of $q$.

***Lemma* 18.** If a function $f : \mathbb{Z} \to \mathbb{C}$ satisfies $f(n) \equiv 0 \bmod q$ for $\gcd(n, q) \neq 1$, is periodic, and completely multiplicative, then $f(n) = \chi(n)$ is a character modulo $q$.

## Properties of Nontrivial Characters
(i) $\chi(1) = 1$ and $\chi(s) \neq 1$, for $\gcd(s, q) > 1$,
(ii) $\chi(st) = \chi(s)\chi(t)$, multiplicative,
(iii) $\chi(s + m) = \chi(s)$, periodic of period $m \geq 1$,
(iv) $|\chi(s)| = 1$, a point in unit circle.

There are several cases depending on the decomposition of the integer $q$.

For $q = p^v$, where $p > 2$ is a prime and $v \geq 1$, each nontrivial character is realized by an equation of the form

$$\chi_s(n) = \begin{cases} e^{i2\pi st/\varphi(q)} & \text{if } \gcd(n,q) = 1, \\ 0 & \text{if } \gcd(n,q) \ne 1, \end{cases} \tag{19}$$

where $t$ = discrete logarithm of $n$ in the multiplicative group of $\mathbb{Z}_q$. Here, $0 \le s, t < \varphi(q)$. The specific case $s = 0$ is called the trivial character, which is defined by

$$\chi_0(n) = \begin{cases} 1 & \text{if } \gcd(n,q) = 1, \\ 0 & \text{if } \gcd(n,q) \ne 1. \end{cases} \tag{20}$$

For $q = 2^v$, $v \ge 1$, a character is realized by one of the three forms described below.
Case $v = 1$, there is a single character in $\hat{G} = \{1 = \chi_0\}$, and it is defined by

$$\chi(n) = \begin{cases} 1 & \text{if } n \equiv 1 \bmod 2, \\ 0 & \text{otherwise.} \end{cases} \tag{21}$$

Case $v = 2$, there are two characters in $\hat{G} = \{1 = \chi_0, \chi\}$, the nontrivial character is defined by

$$\chi(n) = \begin{cases} 1 & \text{if } n \equiv 1 \bmod 4, \\ -1 & \text{if } n \equiv 3 \bmod 4, \\ 0 & \text{if } n \equiv 0, 2 \bmod 4. \end{cases} \tag{22}$$

Case $v > 2$, there are $\varphi(q)$ characters in $\hat{G} = \{1 = \chi_0, \chi_1, ..., \chi_{\varphi(q)-1}\}$, and a nontrivial character is defined by

$$\chi_s(n) = \begin{cases} (-1)^\delta e^{i2\pi st/\varphi(2^{v-3})} & \text{if } n \equiv 1 \bmod 2, \\ 0 & \text{if } n \equiv 0 \bmod 2, \end{cases} \tag{23}$$

for some $0 \le s, t < \varphi(q)$. The integer $n$ is represented as $n \equiv (-1)^\delta 5^t \bmod 2^v$ in the multiplicative group of units $\{-1, 1\} \times \{5^t : 0 \le t < \varphi(2^{v-2})\}$ of $\mathbb{Z}_q$, where $\delta = 0$ if $n \equiv 1 \bmod 4$ or $\delta = 1$ if $n \equiv 3 \bmod 4$. This is due to the fact that this multiplicative group is not cyclic.

A character $\chi$ is *even* if $\chi(t) = \chi(-t)$, otherwise $\chi(t) = -\chi(-t)$, and the character is *odd*. The binary variable $\delta_\chi = 0, 1$ tracks the even odd condition, specifically, $\chi(-n) = (-1)^{\delta_\chi} \chi(n)$. A quadratic character always satisfies $\chi(-1) = -1$. A character $\chi$ is *primitive* if no proper subgroup of the group $\hat{G} = \{1 = \chi_0, \chi_1, ..., \chi_{\varphi(q)-1}\}$ contains it. Under this condition the conductor of a character is the integer $f_\chi = q$.

**Proposition 19.** (Orthogonal relations) For $a \ge 1$, and the set of characters modulo $q$, the followings hold.

$$\text{(i)} \sum_{1 \le n \le q} \chi(n) = \begin{cases} \varphi(q) & \text{if } \chi = \chi_0, \\ 0 & \text{if } \chi \ne \chi_0. \end{cases} \qquad \text{(ii)} \sum_\chi \chi(n)\overline{\chi}(a) = \begin{cases} \varphi(q) & \text{if } n \equiv a \bmod q, \\ 0 & \text{otherwise.} \end{cases} \tag{24}$$

For each character $\chi$ modulo $q$, a Gauss sum is defined by the exponential sum

$$\tau_a(\chi) = \sum_{n=1}^{q} \chi(n) e^{i2\pi an/q}, \quad 0 \le a < q.$$

**Lemma 20.** Let $\chi \ne \chi_0$ be a nontrivial character modulo $q$, and let $\tau(\chi) = \tau_1(\chi)$. Then
(i) $\tau_a(\chi) = \overline{\chi}(a)\tau(\chi)$,  (ii) $\tau_a(\chi)\tau_a(\overline{\chi}) = q$,

(iii) $\tau_a(\chi) = \begin{cases} q^{1/2} & \text{if } q \equiv 1 \bmod 4, \\ iq^{1/2} & \text{if } q \equiv 3 \bmod 4. \end{cases}$

**Lemma 21.** (Polya-Vinogradov) Let $\chi$ be a primitive character modulo $q$, and let $x \le q$. Then

$$\sum_{1 \le n \le x} \chi(n) < \sqrt{q} \log q.$$

These and many other exponential sums arise in the analysis of the zeta and $L$-functions.

**Dirichlet $L$-Functions**

Let $\chi$ be a character modulo $q$. A Dirichlet $L$-function is defined by

$$L(s, \chi) = \sum_{n=1}^{\infty} \frac{\chi(n)}{n^s} = \prod_{\gcd(p,q)=1} \left(1 - \chi(p)/p^s\right)^{-1}, \tag{25}$$

where $s \in \mathbb{C}$, and its inverse is defined by

$$\frac{1}{L(s, \chi)} = \sum_{n=1}^{\infty} \frac{\mu(n)\chi(n)}{n^s} = \prod_{\gcd(p,q)=1} \left(1 + \chi(p)/p^s\right)^{-1}. \tag{26}$$

Both of these power series are absolutely convergent on the half plane $\Re(s) > 1$. Each term $L_p(s, \chi) = \left(1 - \chi(p)p^{-s}\right)^{-1} = \sum_{n \ge 1} \chi(p^n) p^{-ns}$ in the product is a local $L$-function at $p$. The L-function attached to the trivial character has the special form $L(s, \chi_0) = \sum_{n=1}^{\infty} \frac{\chi_0(n)}{n^s} = \zeta(s) \prod_{p \mid q} \left(1 - 1/p^s\right)^{-1}$.

*Example* **22.** Some examples of $L$-functions.

(i) $L(s, \chi_3) = \sum_{n=1}^{\infty} \frac{\chi_3(n)}{n^s}$ for the quadratic field $\mathbb{Q}(\sqrt{-3})$ with the character $\chi_3(n) = \begin{cases} 1 & \text{if } n \equiv 1 \bmod 3, \\ -1 & \text{if } n \equiv 2 \bmod 3, \\ 0 & \text{if } n \equiv 0 \bmod 3. \end{cases}$

(ii) $L(s, \chi_4) = \sum_{n=1}^{\infty} \frac{\chi_4(n)}{n^s}$ for the quadratic field $\mathbb{Q}(\sqrt{-1})$ with the character $\chi_4(n) = \begin{cases} 1 & \text{if } n \equiv 1 \bmod 4, \\ -1 & \text{if } n \equiv 3 \bmod 4, \\ 0 & \text{if } n \equiv 0, 2 \bmod 4. \end{cases}$

(iii) $L(s, \chi_5) = \sum_{n=1}^{\infty} \frac{\chi_5(n)}{n^s}$ for the quadratic field $\mathbb{Q}(\sqrt{-5})$ with the character $\chi_5(n) = \begin{cases} i^a & \text{if } n \equiv a \bmod 5, \\ 0 & \text{if } n \equiv 0 \bmod 5. \end{cases}$

The next result provides some information on the analytic continuation of $L(s, \chi)$ to the complex half plane $\Re e(s) > 0$.

**Theorem 23.** (i) For $\chi_0 = 1$, the $L$-function $L(s,\chi_0)$ has a pole at $s = 1$, and residue $\varphi(q)/q$.

(ii) For $\chi \neq 1$, the representation $L(s, \chi) = \eta_\chi \dfrac{\varphi(q)}{q} \dfrac{s}{s-1} + s \int_1^\infty \dfrac{R(s)}{x^{s+1}} dx$ of $L(s,\chi)$ is holomorphic on the half plane $\Re e(s) > 0$, where $\eta_\chi = 1$ if $\chi(-1) = 1$, otherwise $\eta_\chi = 0$.

The proof is by partial summation applied to $L(s, \chi) = \sum_{n \geq 1} \chi(n) n^{-s}$, see [EL, 227]. Analytic continuation of $L$-functions is a central theme in the theory of modular forms, confer the literature for advanced details.

**Proposition 24.** Let $\chi$ be a character modulo $q$ and let $\Re e(s) > 1$. Then

$$-\frac{L'(s,\chi)}{L(s,\chi)} = \sum_{n=1}^{\infty} \frac{\chi(n)\Lambda(n)}{n^s}. \tag{27}$$

Proof: Start with the Euler product representation $L(s, \chi) = \prod_{\gcd(p,q)=1} \left(1 - \chi(p) p^{-s}\right)^{-1}$, compute the logarithmic derivative $-\dfrac{L'(s,\chi)}{L(s,\chi)} = \sum_{\gcd(p,q)=1} \sum_{n \geq 1} \dfrac{\chi(p)^n \log p}{np^{sn}}$ and simplify it. ∎

**Functional Equation of $L$- Functions**

The functional equation of the Dirichlet $L$-function is given by

$$\begin{aligned}\Lambda(s,\chi) &= q^{(s+\delta)/2} \Gamma((s+\delta)/2) L(s,\chi) \\ &= \varepsilon(\chi) q^{1/2-s} \Lambda(1-s, \overline{\chi}),\end{aligned} \tag{28}$$

where $s \in \mathbb{C}$, the constant

$$\varepsilon(\chi) = i^{-\delta} q^{-1/2} \sum_{t=1}^{q} \chi(t) e^{i 2\pi t/q} \tag{29}$$

(often called the root number) has the absolute value $|\varepsilon(\chi)| = 1$, and $\delta = \delta_\chi = 0$ if $\chi$ is an even character. Otherwise $\delta_\chi = 1$ for odd character. Here it was assumed that $q$ is the conductor of $\chi$. The functional equation can be proved in a similar manner as the functional equation of the zeta function, for example, using the modular theta functions.

**Theorem 25.** Let $L(s, \chi)$ be an $L$-function. Then it has a Weierstrass product expansion

$$s(s-1)^r \Lambda(s,\chi) = e^{A_\chi + B_\chi s} \prod_{\rho \neq 0, 1} (1 - s/\rho) e^{s/\rho}, \tag{30}$$

where $A, B \in \mathbb{C}$ are constants, $r \geq 1$ is the rank, and $\rho$ ranges over the zeros of $L(s, \chi)$ in the critical strip $0 \leq \Re e(s) \leq 1$.

The trivial zeros of the equation $L(s,\chi) = 0$ are $s = -2n - \delta_\chi$, $n \in \mathbb{N}$, and the nontrivial zeros $s = \beta + it$, are on the critical strip $0 \leq \Re e(s) \leq 1$.

**Conjecture 26.** (Generalized Riemann Hypothesis) The nontrivial zeros of the $L$-function are on the line $\Re(\rho) = 1/2$.

**Theorem 27.** (Dirichlet 1837) If $\chi_0 \ne 1$ is a character modulo $q$, then $L(1,\chi) \ne 0$.

In general, it is known that every $L$-function attached to a character $\pi$ satisfies the nonvanishing condition $L(1,\pi) \ne 0$ on the half-plane $\Re(s) \ge 1$.

**Special Values**
The determination of special values of $L(s, \chi)$ in closed forms and the set of critical zeros $\{ \beta + it : L(\beta + it,\chi) = 0 \text{ and } 0 \le \beta \le 1 \}$ are some of the most important open problems in the theory of the $L$-functions.

**Theorem 28.** For a primitive character $\chi$ modulo $q$ and $n \ge 1$ such that $n \equiv (1 - \chi(-1))/2 \mod 2$, the following holds.
(i) $L(1-n, \chi) = -B_{n,\chi}/n$,
(ii) $L(n, \chi) = -(i2\pi q^{-1})^n \tau(\chi) B_{n,\bar{\chi}}/(2n)!$,
where $B_{n,\chi}$ is the $n$th $\chi$-Bernoulli number.

**Zeros Detection.** For all real number $t \in \mathbb{R}$, the real-valued function $Z(t, \chi) = e^{i\theta(t)} L(1/2 + it, \chi)$ has the same zeros lying on the critical line $\Re(s) = 1/2$ as the $L$-function $L(1/2 + it,\chi)$.

A basic rule of zeros detection is the Gram's law. This law claims that there is a zero in the interval $[g_n, g_{n+1}]$, where point $g_n$ satisfies $\theta(g_n) = n\pi$. This is due to the relation $\zeta(1/2 + it, \chi) = Z(t, \chi)\cos\theta(t) + iZ(t, \chi)\sin\theta(t) = (-1)^n Z(t, \chi)$ at the point $t = g_n$. This relation specifies at least one change of sign in the real values of the zeta function. But this law is known to fail infinitely often, consult the literature for finer details.

**Other $L$-Functions**
The zeta function and the Dirichlet $L$-functions are $L$-functions of degree 1 on $GL(1)$. The $L$-functions of degree $N \ge 2$ arise from the automorphic forms on the linear group $GL(N)$ of $N \times N$ matrices. The best known examples of an $L$-functions of degree 2 on $GL(2)$ is

$$L(s,\Delta) = \sum_{n=1}^{\infty} \frac{\tau(n)}{n^{s+11/2}} = \prod_p \left(1 - \tau(p)/p^{s+11/2} + p^{2s}\right)^{-1}, \tag{31}$$

and its functional equation

$$\Lambda(s,\Delta) = \pi^{-s}\Gamma((s+11/2)/2)\Gamma((s+13/2)L(s,\Delta). \tag{32}$$

This $L$-function arises from the discriminant automorphic form

$$\Delta(z) = \sum_{n=1}^{\infty} \tau(n)e^{i2\pi nz} = e^{i2\pi z}\prod_{n\ge 1}\left(1 - e^{i2\pi nz}\right)^{24}, \tag{33}$$

of weight 12 on the modular group $GL(2)$, viz, it satisfies $\Delta((az + b)/(cz + d)) = (cz + d)^{12}\Delta(z)$. Confer the literature for other examples of degree 2 and higher rank $L$-functions.

**3 SURVEY OF THE PRIME NUMBER THEOREM**
The earliest recorded result on the set of prime numbers is a statement on its infinitude.

**Theorem 29.** (Euclid 330 BC)   There are infinitely many prime numbers.

There are about eight variations of the basic time tested combinatorial proof given in [E], and several analytic proofs completely different from the original combinatorial proof, see [RN, p. 3]. The first analytic proof was introduced about three centuries ago.

**Theorem 30.** (Euler 1737)   There are infinitely many prime numbers.

Proof: The verification as given in [EU, p. 240] is something like this:

$$\log \prod_p (1-1/p)^{-1} = \sum_p 1/p + C, \tag{34}$$

where $C > 0$ is a constant. On the other hand,

$$\prod_p (1-1/p)^{-1} = \sum_p 1/n = \log \infty. \tag{35}$$

Therefore $\sum_p 1/p = \log\log\infty$.  ∎

It should be noted that this formal analysis was done over 100 years before Meterns's formula $\sum_{p \leq x} 1/p = \log\log x + c + O(1/\log x)$.

The prime counting function is defined by $\pi(x) = \#\{ p \leq x : p \text{ prime and } x > 1 \}$, and theta and psi functions are defined by

$$\vartheta(x) = \sum_{p \leq x} \log p \quad \text{and} \quad \psi(x) = \sum_{p^n \leq x, n \geq 1} \log p, \tag{36}$$

respectively.

**Theorem 31.** (Chebyshev 1850)   Let $x > 1$. Then
(i) $ax < \vartheta(x) < bx$ for some constants $a, b > 0$,       (ii) $\psi(x) = \vartheta(x) + O(x^{1/2})$,
(iii) $\log x! = \sum_{n \leq x} \psi(x/n)$.

Proof: The first is derived from the properties of the binomial coefficients $(n,k) = n!/k!(n-k)!$, and the second follows from the finite series

$$\psi(x) = \vartheta(x) + \vartheta(x^{1/2}) + \vartheta(x^{1/3}) + \cdots + \vartheta(x^{1/\log x}) = \vartheta(x) + \vartheta(x^{1/2}) + O(\vartheta(x^{1/3})\log x) = \vartheta(x) + O(x^{1/2}). \quad \blacksquare$$

The error terms $\vartheta(x) - x$ and $\psi(x) - x$ are oscillating functions that tend to infinity as $x$ tends to infinity. The peaks and valleys of the oscillations of the latter function are known to satisfy $|\psi(x) - x| > c\sqrt{x}$, $c > 0$ constant, infinitely often. Further, the difference $\psi(x) - \vartheta(x) = O(x^{1/2})$ has a striking resemblance to a *discontinuous* square root function, viz, $\psi(x) - \vartheta(x) \asymp \sqrt{x}$ almost everywhere.

**Theorem 32.** (Littlewood 1914)   The error term satisfies $|\psi(x) - x| = \Omega_{\pm}(x^{1/2}\log\log x)$.

**Theorem 33.**   Let $x \geq 2$, then $\vartheta(x) = x + O(xe^{-(\log x)^{1/7}(\log\log x)^{-2}})$ and $\psi(x) = x + O(xe^{-(\log x)^{1/7}(\log\log x)^{-2}})$.

Over the years a series of improvements on the error term of the psi function has been made, see [IC, p. 39] for more details. Also confer [RS], [SC] and [DT] for other sharper estimates. Some details on the early history of these estimates appears in [NW].

***Theorem* 34.** ([SC])  Assume the Riemann hypothesis, then error terms satisfy
(i) $|\psi(x)-x| < (8\pi)^{-1} x^{1/2} \log^2 x$.
(ii) $|\vartheta(x)-x| < (8\pi)^{-1} x^{1/2} \log^2 x$.

***Theorem* 35.**  Assuming the Riemann hypothesis, the followings hold.
(i) The function $\psi(x)-x$ changes sign in the interval $[x, 19x]$, $x \geq 1$ at least ounce.
(ii) The function $\psi(x)-x$ changes sign in the interval $[x, 2.02x]$, $x \geq 1$ at least ounce.

Confer [MT] for the analysis. A weaker result on the sign changes is given in [NW, p. 226]. It is believed that

$$\liminf_{x \to \infty} \frac{\psi(x)-x}{x^{1/2}(\log\log\log x)^2} = -\frac{1}{2\pi} \quad \text{and} \quad \limsup_{x \to \infty} \frac{\psi(x)-x}{x^{1/2}(\log\log\log x)^2} = \frac{1}{2\pi}.$$

**Prime Number Theorem**
The literature explains that about the same time (circa 1790) both Gauss and Legrendre proposed similar approximations of the form $\pi(x) \approx x/\log x$ for the number of prime up to a given number $x > 1$.

***Conjecture* 36.**  (Gauss 1849)  The number $\pi(x)$ of primes up to $x > 1$ is approximately $\pi(x) \approx \int_2^x (\ln t)^{-1} dt$.

The work of Chebychev established the correct order of magnitude of the function $\pi(x) \asymp x/\log x$, but not the asymptotic form $\pi(x) \sim x/\log x$.

***Theorem* 37.**  (Chebychev 1848)  The number $\pi(x)$ of primes up to $x > 1$ is in the range $\dfrac{ax}{\log x} < \pi(x) < \dfrac{bx}{\log x}$, where $a, b > 0$ are constants.

***Theorem* 38.**  (Prime Number Theorem)  The following statements are equivalents.
(i) $\pi(x) \sim x/\log x$
(ii) $\vartheta(x) \sim x$
(iii) $\psi(x) \sim x$.

The prime number theorem was independently proved by DelaValle Poussin and Hadamard, circa 1896. There are several proofs based on complex analysis, see [NW]. In addition, elementary proofs were independently discovered by both Selberg and Erdos in 1949, see [SG]. The contemporary leading experts did not predict nor expected the existence of an elementary proof.

The Riemann explicit formula

$$\pi(x) = li(x) - \sum_\rho \left(li(x^\rho) + li(x^{1-\rho})\right) + \log \xi(0) + \int_x^\infty \frac{dt}{t(t^2-1)\log t} \tag{37}$$

clearly shows the influence of the zeros of the zeta function on the prime counting function $\pi(x)$. Here the real and complex logarithm integrals are defined by

$$li(x) = \int_2^x \frac{dt}{\ln t} \text{ for } x \in \mathbb{R}, \text{ and } li(e^z) = \int_{-\infty+iv}^z e^t \frac{dt}{t} \text{ for } z = u + iv \in \mathbb{C}, v \neq 0, \tag{38}$$

respectively. The current version of the prime number theorem has a better error term and it is as follows.

***Theorem 39.*** Let $x > 1$. Then $\pi(x) = li(x) + O(xe^{-c(\log x)^{3/5}/(\log \log x)^{-1/5}})$, $c > 0$ constant.

The proof of this version of the prime number theorem is given in [IV, p. 307], the constant $c = .2018$ is given in [FD]. The Riemann hypothesis calls for $\pi(x) = li(x) + O(x^{1/2} \log^2 x)$.

***Theorem 40.*** (Littlewood 1914) (i) The error term $\pi(x) - x/\log x = \Omega(x/(\log x)^2)$ is strictly positive.
(ii) The error term $\pi(x) - li(x) = \Omega_\pm(x^{1/2} \log\log\log x/\log x)$ is oscillatory.

A full proof is given in [EL, p. 191] and a weaker result appears in [MV], and [IV].

***Theorem 41.*** ([SC]) Assuming the Riemann hypothesis, the error term satisfies

$$|\pi(x) - li(x)| < (x^{1/2} \log x)/8\pi. \tag{39}$$

***Theorem 42.*** The followings statements are equivalent.
(i) $\pi(x) \sim x/\log x$,   (ii) $p_n \sim n\log n$.

Proof: To verify that (i) implies (ii), let $p_n$ be the $n$th prime. Then

$$\frac{p_n}{\log p_n} = \pi(p_n) + o(\pi(p_n)) = n + o(n).$$

This in turn implies that $p_n = n\log p_n + o(n\log p_n) = n\log n + o(n\log n)$. Conversely, to verify that (ii) implies (i), take $p_n \leq x < p_{n+1}$. Then the inequalities

$$\frac{p_n}{\log p_n} \leq \frac{x}{\log x} < \frac{p_{n+1}}{\log p_{n+1}}$$

hold. Moreover, since $\pi(x) = n$, the dual inequalities

$$\frac{p_n}{n\log p_n} \leq \frac{x}{\pi(x)\log x} < \frac{p_{n+1}}{n\log p_{n+1}}$$

also hold. Hence, using $p_n \sim n\log n$ and taking the limits yield

$$\lim_{n\to\infty} \frac{p_n}{n\log p_n} = 1 \leq \lim_{n\to\infty} \frac{x}{\pi(x)\log x} < \lim_{n\to\infty} \frac{p_{n+1}}{n\log p_{n+1}} = 1. \blacksquare$$

## 4 REPRESENTATIONS OF ARITHMETIC FUNCTIONS
This section considers the representations of some arithmetic functions. The representations use finite sums, integrals and other related functions. These are elementary techniques but very effective in their realm of applications.

***Theorem 43.*** Let $x > 2$. Then

(i) $\pi(x) = \dfrac{\vartheta(x)}{\log x} + \int_2^x \dfrac{\vartheta(t)}{t \ln^2 t} dt$, (ii) $\vartheta(x) = \pi(x) \log x - \int_2^x \dfrac{\pi(t)}{t \ln^2 t} dt$, (40)

(iii) $\displaystyle\sum_{2 \le n \le x} \dfrac{\Lambda(n)}{\log n} = \dfrac{\psi(x)}{\log x} + \int_2^x \dfrac{\psi(t)}{t \ln^2 t} dt$, (iv) $\displaystyle\sum_{p \le x} 1/p = \dfrac{\pi(x)}{x} + \int_2^x \dfrac{\pi(t)}{t^2} dt$,

(v) $\displaystyle\sum_{2 \le n \le x} \Lambda(n) \log n = \psi(x) \log x - \int_2^x \dfrac{\psi(t)}{t} dt$.

Proof (i): This is done in reverse:

$$\int_2^x \dfrac{\vartheta(t)}{t \ln^2 t} dt = \int_2^x \dfrac{\sum_{p \le t} \ln p}{t \ln^2 t} dt = \sum_{p \le x} \int_p^x \dfrac{\ln p}{t \ln^2 t} dt = -\sum_{p \le x} \dfrac{\ln p}{\ln t}\bigg|_p^x = \pi(x) - \dfrac{\vartheta(x)}{\log x}.$$

Moreover, replacing the $\psi(x) = \vartheta(x) + O(x^{1/2})$ in (i) yields the well known translation formula

$$\pi(x) = \dfrac{\psi(x)}{\log x} + \int_2^x \dfrac{\psi(t)}{t \ln^2 t} dt + O(x^{1/2}/\log x).$$

The corresponding representations over arithmetic progressions are as follow.

**Theorem 44.** Let $x > 2$. Then

(i) $\pi(x, q, a) = \dfrac{\vartheta(x, q, a)}{\log x} + \int_2^x \dfrac{\vartheta(t, q, a)}{t \ln^2 t} dt$, (ii) $\vartheta(x, q, a) = \pi(x, q, a) \log x - \int_2^x \dfrac{\pi(t, q, a)}{t \ln^2 t} dt$, (41)

(iii) $\displaystyle\sum_{p \le x} \dfrac{\chi(p)}{p} = \dfrac{\vartheta(x, \chi)}{x \log x} + \int_2^x \vartheta(t, \chi) \dfrac{1 + \ln t}{t^2 \ln^2 t} dt$, (iv) $\displaystyle\sum_{p \le x, p \equiv a \bmod q} 1/p = \dfrac{\pi(x, q, a)}{x} + \int_2^x \dfrac{\pi(t, q, a)}{t^2} dt$.

The proofs of these identities are based on partial summations. As an illustration, the verification of (ii) is as follows: For a character $\chi$ modulo $q$, define the $\chi$-twist of the theta function as $\vartheta(x, \chi) = \sum_{p \le x} \chi(p) \log p$, and put $S(x) = \vartheta(x, \chi)$ and $g(n) = n \log n$. Now using the basic partial summation formula

$$\sum_{n \le x} f(n) g(n) = S(x) g(x) - \int_x^\infty S(t) dg(t),$$ (42)

where $S(x) = \sum_{n \le x} f(n)$, see [AP, p. 78], [HW, p. 460] for the basic technique. Lastly, one has

$$\sum_{p \le x} \dfrac{\chi(p)}{p} = \sum_{p \le x} \dfrac{\chi(p) \log p}{p \log p} = \dfrac{\vartheta(x, \chi)}{x \log x} - \int_2^x \vartheta(t, \chi) \dfrac{d}{dt}\left(\dfrac{1}{t \ln t}\right).$$ (43)

**Theorem 45.** (i) Let $x \ge 1$, then $\sum_{p \le x} 1/p = \log \log x + B + O(1/\log^2 x)$.

(ii) Assume the Riemann hypothesis, then $\sum_{p \le x} 1/p = \log \log x + B + O(\log x / x^{1/2})$,

where $B = .2614972128\ldots$, is Mertens constant.

The first is due to Mertens, 1874, and the second is due to the combine works of several authors, refer to the literature.

# 5 PRIMES IN SHORT INTERVALS

This section begins with a survey of the analysis of primes in short intervals and ends with a proposed improvement of this analysis. There is a vast and extensive literature on these topics, and an effort was made to provide sufficient references to the literature. The leading techniques are the zero density methods and the sieve methods, see [HB] and [PT] for surveys.

**Tchebychev Method**

The earliest result on primes in short intervals was achieved by Tchebychev, it deals with intervals of the form $[x, x + y]$ with $y = cx$, $c \approx 2$. The analysis was further developed to $y \asymp x$, see [HW, p. 455], [RN]. The special case of $c = 2$ is widely known as Bertrand's Postulate, and the case $y \asymp x$ follows from the prime number theorem.

***Proposition* 46.** (Tchebychev 1852) The interval $[x, 2x]$ contains a prime for all $x \geq 1$.

Proof: The analysis springs from the two simple properties of the central binomial coefficient $(2n, n)$:
(i) $(2n, n) = (2n)!/(n!)^2$ is not divisible by any prime $p$ in the range $2n/3 \leq p \leq n$, and
(ii) $(2n, n) > 2^{2n}/(2n+1)$ is larger than the mean value of the binomial coefficients.
Firstly, assume that the interval $[n, 2n]$ contains no prime, then

$$\binom{2n}{n} \leq \prod_{p^v \leq 2n, v \geq 1} p^v \leq \left(\prod_{p \leq 2n/3} p\right)\left(\prod_{p^v \leq 2n, v \geq 2} p^v\right) \leq 4^{2n/3 + O(n/\log n)} \times (2n)^{\sqrt{2n}}, \tag{44}$$

this follows from properties (i), $\prod_{p \leq x} p \leq e^{\vartheta(x)} \leq e^{x + O(x/\log x)}$, and the fact that there are less than $\sqrt{2n}$ prime powers $p^v$, $v \geq 2$, in the interval $[n, 2n]$. Lastly, by property (ii) one has

$$2^{2n}/(2n+1) < \binom{2n}{n} \leq 4^{2n/3} \times (2n)^{\sqrt{2n}}. \tag{45}$$

Clearly this last inequality is a contradiction for large $n$, the small integers $n$ cases are resolved by hand calculations. ∎

It appears that the same binomial coefficient argument can not be extended to smaller intervals $[n, n + k]$ with small $k < n$ since the inequality $2^{2n}/(2n+1) < (2n)!/(k!(n-k)!)$ for the noncentral binomial coefficients $(2n, k) = (2n)!/(k!(n-k)!)$ is no longer true for small $k = n^\delta$, $\delta > 0$.

There are other ways of proving Proposition 46. For example, a single line proof is also possible, and it is as follows: Since $ax < \vartheta(x) < bx$ with $2a > b$, see Proposition 31, the difference $\vartheta(x + y) - \vartheta(x) > (2a - b)x > \log(x)$ implies that the interval $(x, x + y]$ contains at least one prime.

***Proposition* 47.** Let $\mu > 0$ be a small fixed number, and let $y = \mu x$. Then the interval $[x, x + y]$ contains a prime for all sufficiently large $x \geq 1$.

Proof: Using a weak version of the prime number theorem, one has

$$\pi(x+y) - \pi(x) = \frac{x+y}{\log(x+y)} + O(\frac{x+y}{\log^2(x+y)}) - \left(\frac{x}{\log x} + O(\frac{x}{\log^2 x})\right) \quad (46)$$

$$= \frac{\mu x}{\log(x+y)}(1+o(1)) \geq 1.$$

∎

Partial summation techniques can be used to improve the Tchebychev result to smaller intervals. The basic idea is illustrated below.

**Theorem 48.** Let $\varepsilon > 0$ be fixed small number and let $x > 0$ be a sufficiently large number. Then the interval $[x, x+y]$ contains a primes for any number $y \geq x/(\log x)^{1-\varepsilon}$.

Proof: By Theorem 45, one has the differencing

$$\sum_{x < p \leq x+y} 1/p \geq \log\log(x+y) + O(1/\log^2(x+y)) - \left(\log\log(x) + O(1/\log^2 x)\right)$$

$$\geq \log(1 + \log(1 + y/x)/\log x) + O(1/\log^2 x) \quad (47)$$

$$\geq \frac{cy}{x \log x} + O(1/\log^2 x),$$

where the constant $c \approx 1$ arises by repeated use of the identity $\log(1+z) \geq c_0 z$ for any complex number $z$ such that $|z| < 1$, and some constant $0 < c_0 < 1$. By setting $y = x/(\log x)^{1-\varepsilon}$, the claim immediately follows from the inequalities

$$\sum_{x < p \leq x+y} 1/p \geq c/\log^{2-\varepsilon} x + O(1/\log^2 x) \geq 1/x.$$

∎

**DelaValle Poussin Method**
The earliest result on primes in short intervals $[x, x+y]$ of subexponential sizes was achieved by DelaValle Poussin. The analysis consists of an improvement in the error term in the prime number theorem. First of all, the asymptotic formula $\pi(x) \sim x/\log x$ implies the existence of primes in the intervals $[x, x+y]$ with $y \asymp x/\log x$. Secondly, the improved error term in $\pi(x) = li(x) + O(xe^{-c\sqrt{\log x}})$ or better, see Theorem 39, implies the existence of primes in the intervals $[x, x+y]$ of subexponential sizes $y = xe^{-c(\log x)^{1/2-\varepsilon}}$ with $c, \varepsilon > 0$ constants.

**Hoheisel Method**
The earliest result on short intervals $(x, x+y]$ of size $y = x^{1-\varepsilon}$, $1/3300 < \varepsilon < 1$, or better was achieved by Hoheisel in 1930, see [KA, p. 100], [PT]. The current record is approximately $y = x^{1/2+\varepsilon}$, $1/21 < \varepsilon < 1$. The analysis, often called the *zero density method*, revolves around the explicit formula in Theorem 49. Since the zeros of the zeta function are involved in the explicit formula, the analysis is difficult and requires extensive knowledge of the zero-free region, zero-density estimates and so on. A few prominent results are given here.

**Theorem 49.** (Landau 1906) Let $\{\rho = \beta + i\gamma\}$ be the critical zeros of $\zeta(s)$. Then

$$\psi(x) = x - \sum_{|\rho| \leq T} \frac{x^\rho}{\rho} + O(xT^{-1}(\log Tx)^2) + O(\log x), \quad (48)$$

uniformly for $T > T_0$.

This is an effective form of Theorem 9, for a proof, see [IV, p. 301], and [DP]. Any result that proves $\psi(x+y) - \psi(x) = (1+o(1))y > 0$ for $0 < y < x$, immediately implies the existence of primes in the short intervals $[x, x+y]$, see [IV, 315], [IC, p. 264], and [KA, p. 100] for discussions. Typical results for short the intervals are as shown below. These results illustrate the basic technique well at an elementary level.

**Theorem 50.** Let $y \geq x^{3/4} e^{\log^8 x}$, $x > x_0$. Then the asymptotic formula $\psi(x+y) - \psi(x) = y + O(ye^{-\log^{-1} x})$ holds.

*Sketch of the proof*: Put $\rho = \beta + it$, $0 \leq \beta \leq 1$, and let $t \in \mathbb{R}$. Since

$$\left| \frac{(x+y)^\rho - x^\rho}{\rho} \right| = \left| \int_x^{x+y} t^{\rho-1} dt \right| \leq \int_x^{x+y} t^{\beta-1} dt \leq x^{\beta-1} y, \tag{49}$$

the sum of these terms is

$$\sum_{|\rho| \leq T} x^\beta = \sum_{|\rho| \leq T} \log x \int_0^\beta x^t dt + 1 = N(T) + \log x \sum_{|\rho| \leq T} \int_0^\beta x^t dt, \tag{50}$$

where $N(T) = O(T \log T)$ is the number of zero in the rectangle $|\rho| \leq T$. This is followed by a not too difficult density theorem step to reduce the last finite sum/integral to subexponential size in $x$, that is, $\sum_{|\rho| \leq T} x^\beta = O(xe^{-\log^{-1} x})$, see [KR, p. 98] for the complete details. ∎

This result immediately implies the existence of primes in the interval $[x - x^{3/4+\varepsilon}, x]$. The best result using zero density method appears to be the following.

**Theorem 51.** ([HL]) Let $\theta > 7/12$. Then the following holds.
(i) The interval $[x, x+y]$ contains $y(1 + o(1))/\log x$ primes for $y \geq x^\theta$.
(ii) Let $p_n$ be the $n$th prime, then $p_{n+1} - p_n \leq p_n^\theta$.

Proof, confer [HN, p. 120]: For any zero $\rho = \beta + it$, $0 \leq \beta \leq 1$, $t \in \mathbb{R}$, and $y \leq x$, one has

$$\frac{(x+y)^\rho - x^\rho}{\rho} = \int_x^{x+y} t^{\rho-1} dt. \tag{51}$$

From this it follows that

$$\left| \sum_{|\rho| \leq T} \frac{(x+y)^\rho - x^\rho}{\rho} \right| \leq 2x^{\beta-1} y N(\sigma, T), \tag{52}$$

where the index ranges over the zeros with $\beta \geq \sigma$, and $|\rho| \leq T$. Now put $T = x^{1-\theta} \log^3 x$ in the density estimate $N(\sigma, T) = T^{12(1-\sigma)/5} \log^B T$, $B > 0$ constant, see Theorem 15. In a zero-free region

$$\sigma > \sigma_0 = 1 - \frac{c}{(\log T)^{2/3} (\log \log T)^{1/3}}, \tag{53}$$

where $c > 0$ is a constant, the zero count $N(\sigma, T) = 0$, and the penultimate finite sum (53) vanishes. But for $\sigma < \sigma_0$, finite sum (53) is

$$\left| \sum_{|\rho| \leq T} \frac{(x+y)^\rho - x^\rho}{\rho} \right| \leq 2x^{\beta-1} y N(\sigma, T)$$

$$\leq 2x^{\sigma-1} y \left( x^{1-\theta} \log^3 x \right)^{12(1-\sigma)/5} \log^{B+1} x \quad (54)$$

$$\leq 2 y x^{12(1-\sigma)(1-\theta)/5 - 1} (\log x)^{36(1-\theta)/5 + B + 1} \leq y / \log x$$

for $\theta > 12/5$ and some $B > 0$. Hence, for $y \geq x^\theta$, one arrives at the equation $\psi(x+y) - \psi(x) = y(1 + o(1))$. ∎

### Sieve Methods

The sieve methods are based on the modern theory of the sieve of Eratothenes, circa 200 BC. These techniques seek effective estimates of the cardinalities of constrained sets of integers. A typical case is the sets of integers

$$S(x, y, z) = \#\{ n : x \leq n \leq x + y, \text{ and } \gcd(n, p) = 1 \text{ for } p \leq z \} \quad (55)$$

and other complicated sets of integers. The most recent result on the application of sieve methods to the theory of primes in short intervals is the following.

**Theorem 52.** ([BK]) For all large $x$, the interval $[x - x^{.525}, x]$ contains primes numbers.

### Conditional Methods for the Short Intervals $[x, x + x^{1/2 + \varepsilon}]$

The Riemann hypothesis limits the size of the short intervals that can be analyzed using the zero density method. This is probably true for the sieve methods too. The best possible result under this hypothesis is the claim that any short interval $[x, x + x^{1/2 + \varepsilon}]$ contains primes.

**Theorem 53.** (von Koch 1901) If the Riemann hypothesis holds, then
(i) There are at least $cx^{1/2} \log x$ primes in the interval $[x, x + x^{1/2 + \varepsilon}]$, $c > 0$ constant.
(ii) The prime gap is of order $p_{n+1} - p_n = O(p_n^{1/2} \log^2 p_n)$.
(iii) $\psi(x + y) - \psi(x) = y + O(x^{1/2} \log^2 x)$.

Proof. These are derived from the integral $\pi(x+y) - \pi(x) = \int_x^{x+y} (\log t)^{-1} dt + O(x^{1/2} \log^2 x)$, where $y = x^{1/2 + \varepsilon}$, see [NW, p. 245], and for a different proof of (ii), see [IV, p. 321]. ∎

The main obstacle in deriving unconditional results appears to be the determination of effective estimates of the exponential sums

$$E(x) = \sum_{|\rho| \leq T} \frac{x^\rho}{\rho} \quad \text{and} \quad E(x, y) = \sum_{|\rho| \leq T} \frac{(x+y)^\rho - x^\rho}{\rho}. \quad (56)$$

For $\theta = \max \{ \Re(\rho) : \zeta(\rho) = 0 \}$, the estimates

(i) $| E(x) | \leq O(x^\theta \log^2 x)$,  (ii) $| E(x, y) | \leq O(x^\theta \log^2 x)$,

can be determined using Theorem 12, which is about the number of zeros $|\rho| \leq T$ with $T^{-1} x \leq x^\theta$, or Theorem 49.

**Lemma 54.** Assume the Riemann hypothesis, then $\sum_{x < p \leq x+y} \log p = y + O(x^{1/2} \log^2 x)$.

To arrive at this result, apply Theorem 43-ii, it also can be derived from Theorem 49. Currently the best unconditional result is

$$\sum_{x<p\leq x+y} \log p = y + O(xe^{-2.18(\log x)^{3/5}(\log\log x)^{1/5}}),\qquad(57)$$

This follows from the current version of the prime number theorem, see Theorem 39.

**Psi-Theta Method**

Partial summation Methods, Sieve Methods and explicit formulas are the standard methods of estimating finite sums of functions over the primes such as

$$S(x) = \sum_{p\leq x} f(p) \quad\text{or}\quad S_\Lambda(x) = \sum_{p\leq x} \Lambda(n) f(n).\qquad(58)$$

A different approach for estimating such finite sums based on the elementary analysis and the prime number theorem will be considered here. This technique is suitable for estimating finite sums of certain class of functions over primes. The result will be used to derive an unconditional result akin to the previous Lemma. It has an interesting application to the theory of primes in short intervals.

For a well constrained real-valued function $f(t) \geq 0$ of bounded variation $0 \leq f(t) \leq M$ on the interval $[\alpha, \beta]$, there exists a finite sum lower estimate of the form

$$\sum_{\alpha<n\leq\beta} f(n) \geq (\beta-\alpha)\min_{\alpha<t\leq\beta} f(t).\qquad(59)$$

Similarly, there exists a finite sum upper estimate of the form $\sum_{\alpha<n\leq\beta} f(n) \leq (\beta-\alpha)\max_{\alpha<t\leq\beta} f(t)$. In some cases the integral approximations $\sum_{\alpha<n\leq\beta} f(n) \geq c_0\int_\alpha^\beta f(t)dt$ and $\sum_{\alpha<n\leq\beta} f(n) \leq c_1\int_a^b f(t)dt$ provide sharper estimates.

**Lemma 55.** Let $x, y$ be sufficiently large real numbers such that $0 \leq y \leq x$. Then
(i) $\sum_{x<p\leq x+y} \log p \geq y + O(y\log\log x/\log x)$.
(ii) $\sum_{x<p\leq x+y} \log p \leq y(1 + 2\log\log x/\log x)$.

*Proof* (i): Let $p_n$ is the $n$th prime, and use the asymptotic expression $p_n \sim n\log n$, see Theorem 42, to obtain the lower estimate

$$\sum_{x<p\leq x+y}\log p \geq \sum_{\alpha<n\leq\beta}\log p_n \geq \sum_{\alpha<n\leq\beta}\log(n\log n + o(n\log n)) \geq \sum_{\alpha<n\leq\beta}\log(n),\qquad(60)$$

where $p_n = n\log n + o(n\log n)$. The approximate limits $\alpha$ and $\beta$ are extracted from the inequalities

$$x/(1 + y/x)) < x < n\log n + o(n\log n) \leq x + y.\qquad(61)$$

As functions of $x$ and $y$, the approximate limits are given by the inequalities

$$\frac{x}{\log n + o(\log n)} \leq \frac{x}{\log x - y/x} < n \leq \frac{x+y}{\log x + y/x} \leq \frac{x+y}{\log n + o(\log n)}.\qquad(62)$$

More precisely, $\alpha(x) = x/(\log x - y/x)$ and $\beta(x) = (x+y)/(\log x + y/x)$, which specifies a slightly smaller range than the actual range of $n$. Continuing the calculation of the lower estimate yields

$$\sum_{\alpha < n \leq \beta} \log(n) \geq (\beta - \alpha) \log \alpha$$

$$= \left( \frac{x+y}{\log x + y/x} - \frac{x}{\log x - y/x} \right) \log \left( \frac{x}{\log x - y/x} \right) \tag{63}$$

$$\geq y - c_0 \frac{y \log \log x}{\log x},$$

where $c_0 > 0$ is a constant. The proof of (ii) is almost the same mutatis mutandis. ∎

A similar proof can be arranged using the inequalities $an \log n \leq p_n \leq bn \log n$ with $a, b > 0$ instead. A light numerical test was conducted on a machine to test the algorithm on random $x < 10^9$ with $y = x^{1/2}$, using $\alpha = x/(\log x - \delta(x))$ and $\beta = (x+y)/(\log x + y/x)$ where $\delta(x) = c_1 x^{-1/2}$, $c_1 > 0$ constant. The result of the test was as predicted.

Lemma 55 provides a rough estimate of the variance of the random variable $\vartheta(x+y) - \vartheta(x) - y$, which is equivalent to a rough estimate of the variance of the random variable $\psi(x+y) - \psi(x) - y$ up to the square root error term. The variance of the latter random variable is believed to be approximately $y \log x$. An abstract of the precise statement as given in [MJ], is as follows.

**Conjecture 56.** Let $\varepsilon > 0$, $N > 1$, and $N^\varepsilon < y < N^{1-\varepsilon}$. For $y \leq x < N$, the random variable $\psi(x+y) - \psi(x) - y$ is normally distributed with mean $\sim y$ and variance $\int_2^N (\psi(x+y) - \psi(x) - y)^2 dx \sim Ny \log N/y$.

The *psi-theta method* arises from remarkable identity $\psi(x) = \vartheta(x) + O(x^{1/2})$ between the psi and the theta functions, see Theorem 31-ii. This identity has the proper square root error term built-in as required by the Riemann hypothesis, for a current perspective on the error terms of various quantities in number theory, confer [MZ]. The main result slightly develops this idea to establish the existence of primes in short intervals.

**Theorem 1.** For all sufficiently large numbers $x > 0$, and $y \geq x^{1/2}(\log x)^\varepsilon$, $\varepsilon > 0$, the interval $[x, x+y]$ contains some primes.

*Proof*: The psi - theta functions identity leads to the differencing

$$\psi(x+y) - \psi(x) = \sum_{x < p \leq x+y} \log p + O(x^{1/2}). \tag{64}$$

Applying Lemma 55 and setting $y = x^{1/2}(\log x)^\varepsilon$, return

$$\psi(x+y) - \psi(x) \geq y + O(y \log \log x / \log x) + O(x^{1/2})$$
$$\geq x^{1/2}(\log x)^\varepsilon \left( 1 - c_2 \frac{\log \log x}{\log x} + c_3 \frac{1}{(\log x)^\varepsilon} \right), \tag{65}$$

where $c_2, c_3 > 0$ are constants. Ergo, for sufficiently large $x$, the interval $[x, x+y]$ contains a prime. ∎

In comparison, the best results in the literature are Theorem 51 using zero density methods, Theorem 52 using state of the art sieve methods, and the almost identical Theorem 53 which assumed the Riemann hypothesis.

# 6 PRIMES IN ARITHMETIC PROGRESSIONS

The prime counting function in an arithmetic progression $\{ qn + a : n \in \mathbb{N} \text{ and } \gcd(a, q) = 1 \}$ is defined by

$$\pi(x, q, a) = \#\{ p \leq x : \text{prime } p \equiv a \mod q \}. \tag{66}$$

The corresponding theta and psi functions are defined by

$$\vartheta(x, q, a) = \sum_{p \leq x, \, p \equiv a \bmod q} \log p \quad \text{and} \quad \psi(x, q, a) = \sum_{p^n \leq x, \, p \equiv a \bmod q} \log p, \tag{67}$$

respectively.

Let $\chi$ be a character modulo $q$. The $\chi$-twists of the theta and psi functions are defined by

$$\vartheta(x, \chi) = \sum_{p \leq x} \chi(p) \log p \quad \text{and} \quad \psi(x, \chi) = \sum_{p \leq x} \chi(n) \Lambda(n) \tag{68}$$

respectively. The orthogonal relation

$$\varphi(q)^{-1} \sum_{\chi} \overline{\chi}(a) \chi(n) = \begin{cases} 1 & \text{if } a \equiv n \bmod q, \\ 0 & \text{otherwise,} \end{cases} \tag{69}$$

of the characters modulo $q$ can be used to link the discrete domain of congruences to the analytic domain of character sums and functions. A simple calculation using the orthogonal relation gives

$$\vartheta(x, q, a) = \varphi(q)^{-1} \sum_{\chi} \overline{\chi}(a) \vartheta(x, \chi) \quad \text{and} \quad \psi(x, q, a) = \varphi(q)^{-1} \sum_{\chi} \overline{\chi}(a) \psi(x, \chi). \tag{70}$$

These character sums are used in the analysis of primes in arithmetic progressions.

A proof of the infinitude of primes in a single arithmetic progression $\{ qn + a : n \in \mathbb{N} \text{ and } \gcd(a, q) = 1 \}$ for a some fixed pair $a$ and $q$ is not to difficult to realize, using the same argument as in Euclid Theorem mutatis mutandis. For example, it is easy to verify that there are infinitely many primes of the forms $3n \pm 2$, $4n \pm 1$, $n \in \mathbb{N}$. The general idea of the possible cases is stated below and a proof appears in [RM].

***Theorem* 57.** (Schur 1912) There exist an Euclidean-type proof for the infinitude of primes in an arithmetic progression $\{ qn + a : n \in \mathbb{N} \text{ and } \gcd(a, q) = 1 \}$ if and only if $a^2 \equiv 1 \mod q$.

For example, this method can not be used to that there are infinitely many primes $p \equiv 2 \mod 5$ since $2^2 \equiv 4 \mod 5$. A much more elaborate complex analysis of the characters and the corresponding $L$-functions leads to the following general result for any pair $a$ and $q$.

***Theorem* 58.** (Dirichlet 1837) The arithmetic progression $\{ qn + a : n \in \mathbb{N} \text{ and } \gcd(a, q) = 1 \}$ has infinitely many primes $p \equiv a \mod q$.

***Theorem* 59.** (*Prime number theorem in arithmetic progressions*) Let $a$ and $q$ be relatively prime integers. Then
(i) $\pi(x, q, a) \sim x / \varphi(q) \log x$,     (ii) $\vartheta(x, q, a) \sim x / \varphi(q)$,     (iii) $\psi(x, q, a) \sim x / \varphi(q)$.

This result is due to delaValle Poussin, a sharper version is as follows.

***Theorem* 60.** (Siegel-Walfisz) For relatively prime integers $a$ and $q$ and a constant $c > 0$, one has

$$\pi(x, a, q) = \frac{1}{\varphi(q)} li(x) + O(xe^{-c(\log x)^{1/2}}) \tag{71}$$

for any modulo $q \leq (\log x)^B$ and any $B > 0$.

For large value of $q > x^\varepsilon$ with any small number $\varepsilon > 0$, there is a much weaker result.

***Theorem* 61.** (Brun-Tichmarch) Let $a$ and $q$ be relatively prime integers. Then

$$\pi(x, a, q) < \frac{2x}{\varphi(q) \log(x/q)} \tag{72}$$

for $x > q$.

The next result is an unconditional result on the error term for the prime number theorem on arithmetic progression, however, it only an on average result. It is sort of a generalized Riemann hypothesis on average.

***Theorem* 62.** (Bombieri-Vinogradov 1966) For $Q \leq x^{1/2 - \varepsilon}$, $\varepsilon > 0$, and $B > 0$,

$$\sum_{q \leq Q, t \leq x} \max_{\gcd(a,q)=1} \left| \psi(t, q, a) - \frac{t}{\varphi(q)} \right| = O\left( \frac{x}{(\log x)^B} + Qx^{1/2} (\log Qx)^4 \right) \tag{73}$$

uniformly for $Q \geq 1$.

A proof is given in [CU, p. 170].

***Theorem* 63.** (Barban et al) For any $A > 0$ and $Q$ such that $x/\log x \leq Q \leq x$, the "variance" of $\psi(x, a, q)$ satisfies the inequality

$$\sum_{q \leq Q} \max_{\gcd(a,q)=1} \left| \psi(x, q, a) - \frac{x}{\varphi(q)} \right|^2 \leq Qx \log x \tag{74}$$

for $x \geq 1$.

**Explicit Formula**
***Theorem* 64.** Let $\{ \rho = \beta + it \}$ be the zeros of $L(s, \chi)$ on the critical strip $0 \leq \Re(s) \leq 1$ and let and $0 < T < x$. Then

$$\psi(x, q, a) = \frac{x}{\varphi(q)} - \frac{1}{\varphi(q)} \sum_\chi \overline{\chi}(a) \sum_{|\rho| \leq T} \frac{x^\rho - 1}{\rho} + O(xT^{-1} \log^2 qx) \tag{75}$$

for any pair of integers such that $\gcd(a, q) = 1$, and $1 \leq T \leq x$.

A proof of this appears in [IC, 121].

**Short Intervals And Arithmetic Progressions**
This last subsection consider a result for primes from arithmetic progressions in short intervals.

**Theorem 2** Let $q \leq 2(\log x)^A$, $A > 0$, and let $x > 1$ be a sufficiently large number. Then the interval $[x, 2x]$ contains a prime number from the arithmetic progression $p \equiv a \bmod q$, $\gcd(a, q) = 1$.

Proof: Given $q = 2(\log x)^A$, $A > 0$, let $B = A + 1$. By Bombieri-Vinogradov Theorem 62 one has

$$\left| \psi(x,q,a) - \frac{x}{\varphi(q)} \right| \leq \frac{c_1 x}{(\log x)^B} \tag{76}$$

uniformly on $a$ and $q$, where $c_1 > 0$ is a constant. Hence, for $y = 2x$ the psi difference resolves to

$$\begin{aligned}\psi(x+y,q,a) - \psi(x,q,a) &= \frac{y}{\varphi(q)} + O\left(x/(\log x)^B\right) \\ &\geq \frac{x}{(\log x)^A}\left(1 + O(\frac{1}{\log x})\right)\end{aligned} \tag{77}$$

The positivity condition follows from the hypothesis $y = 2x$ and $\varphi(q) < q = (\log x)^A < (\log x)^{A+1} = (\log x)^B$. ∎

A recent result in [KR] for the intervals $[e^x, e^{x+\varepsilon}]$ is similar to the above result but the analysis is more complex.

# 7 USEFUL FORMULAS
**Power Series Expansions**
**Logarithm function** $\log(1+x) = x - \dfrac{x^2}{2} + \dfrac{x^3}{3} - \cdots + (-1)^{n+1}\dfrac{x^n}{n} + \cdots$

**Arithmetic Functions**
**Mobius Function**

$$\mu(N) = \begin{cases} (-1)^k & \text{if } N = p_1 p_2 \cdots p_k, \\ 0 & \text{if } p^2 \mid N \text{ for some prime } p. \end{cases}$$

**Mobius Inversion Formulae**
(i) $f(n) = \sum_{d \mid n} g(d) \quad \Leftrightarrow \quad g(n) = \sum_{d \mid n} \mu(d) f(n/d)$,

(ii) $f(n) = \prod_{d \mid n} g(d) \quad \Leftrightarrow \quad g(n) = \prod_{d \mid n} f(d)^{\mu(n/d)}$.

**Euler Function**
(i) $\varphi(N) = N \prod_{p \mid N}(1 - 1/p)$    (ii) $\varphi(n) = N \sum_{d \mid N} \mu(d) \dfrac{1}{d}$.

**von Mangoldt Function**

$$\Lambda(n) = \begin{cases} \log p & \text{if } n = p^m, m \geq 1, \\ 0 & \text{otherwise.} \end{cases}$$

There are many other representations of this function, for example,

$$\Lambda(n) = \sum_{d \mid n} \mu(d) \log n/d = -\sum_{d \mid n} \mu(d) \log d.$$

The extended version is of the form
$$\Lambda_k(n) = \sum_{d \mid n} \mu(d)(\log n/d)^k \text{ for } \omega(n) \le k; \text{ otherwise } \Lambda_k(n) = 0 \text{ for } \omega(n) > k.$$

**Bernoulli Polynomials and Numbers.** The generating function of the Bernoulli numbers is
$$\frac{z}{e^z - 1} = \sum_{n=0}^{\infty} B_n \frac{z^n}{n!}.$$

There are many recursive formulas for computing the 2nth Bernoulli number. Two of these recursive formulas are given by

(i) $\sum_{k=0}^{2n-1} \binom{n}{k} B_k = 0$
(ii) $B_{2n} = \frac{-1}{(n+1)(2n+1)} \sum_{k=1}^{n/2} (2n - 2k + 1) \binom{n+1}{2k+1} B_{2n-2k}$

for $n \ge 1$. The individual number can be approximated using the Euler formula

(i) $B_{2n} = (-1)^{n-1} \frac{2(2n)!}{(2\pi)^n} \zeta(2n)$,
(ii) $|B_{2n}| \sim 2\sqrt{4\pi n} \left(\frac{n}{e\pi}\right)^{2n}$.

A sample of these numbers is listed below.
$B_0 = 1$         $B_1 = -1/2$         $B_2 = 1/6$         $B_4 = -1/30$
$B_6 = 1/42$     $B_8 = -1/30$        $B_{10} = 5/66$     $B_{12} = -691/2730$
$B_{14} = 7/6$   $B_{16} = -3617/510$, $B_{18} = 43867/798$, $B_{20} = -17460/330$, ... .

**Bernoulli Polynomials.** The generating function of the Bernoulli polynomials is
$$\frac{ze^{xz}}{e^z - 1} = \sum_{n=0}^{\infty} B_n(x) \frac{z^n}{n!}.$$

Unlike the nth number $B_n$, the nth polynomial $B_n(x)$ does not vanish at odd parameter $n > 1$. The first few are listed here.
$B_0(x) = B_0$,                              $B_1(x) = x - B_1$,
$B_2(x) = x^2 - x + B_2$,                    $B_3(x) = x^3 - (3/2)x^2 + (1/2)x$,
$B_4(x) = x^4 - 2x^3 + x^2 - B_4$,           $B_5(x) = x^5 - (5/2)x^4 + (5/3)x^3 - (1/6)x$,
$B_6(x) = x^6 - 3x^5 + (5/2)x^4 - (1/2)x^2$, ...,

The recursion $B_n(x) = \sum_{k=0}^{n} \binom{n}{k} B_k x^{n-k}$ gives the nth polynomial.

**Properties**
(i) $B_n'(x) = nB_{n-1}(x)$,
(ii) $B_{n+1}(x+1) - B_{n-1}(x) = nx^{n-1}$,
(iii) $\int_0^1 B_n(x)dx = 0$.

The Fourier series expansion of the nth Bernoulli polynomial on the interval [0,1] is given by
is $B_n(x) = -n! \sum_{-\infty}^{\infty} \frac{e^{i2\pi dx}}{(i2\pi d)^n}$.

In particular,
$$B_{2n}(x) = \frac{(-1)^{n-1} 2(2n)!}{(2\pi)^{2n}} \sum_{d=1}^{\infty} \frac{\cos 2\pi dx}{d^{2n}} \text{ and } B_{2n-1}(x) = \frac{(-1)^n 2(2n-1)!}{(2\pi)^{2n-1}} \sum_{d=1}^{\infty} \frac{\cos 2\pi dx}{d^{2n-1}} \text{ for } 0 \le x \le 1.$$

**χ-Bernoulli Polynomials.** The generating function of the χ-Bernoulli polynomials is

$$\sum_{0 \leq t < q} \chi(t) \frac{ze^{tz}}{e^{qz}-1} = \sum_{n=0}^{\infty} B_{n,\chi}(x) \frac{z^n}{n!}.$$

The constant term $B_{n,\chi} = B_{n,\chi}(0)$ of the $n$th polynomial $B_{n,\chi}(x)$ is called the $n$th χ- Bernoulli number. This complex number coincides with the exponential sum $B_{n,\chi} = q^{n-1} \sum_{0 \leq t < q} \chi(t) B_{n,\chi}(t/q)$.

**Euler-Maclaurin Summation Formula**

Let $f : [a, b] \to \mathbb{C}$ be continuous and differentiable function of order $n$th. Then

$$\sum_{a+1}^{b} f(n) = \int_a^b f(x)\,dx + \sum_{k=1}^{n} (-1)^k \frac{B_k}{k!} \left( f^{(k-1)}(b) - f^{(k-1)}(a) \right) - \frac{(-1)^n}{n!} \int_a^b B_n(\{x\}) f^{(n)}(x)\,dx,$$

where the last term is the error term.

**Abel Summation Formula.** This is an elementary but extremely useful summation technique. Let $f, g : \mathbb{R} \to \mathbb{C}$ be continuously differentiable complex-valued functions on the real line. Let $S(x) = \sum_{n \leq x} f(n)$. Then

$$\sum_{n \leq x} f(n)g(n) = S(x)g(x) - \int_x^{\infty} S(t)\,dg(t). \tag{78}$$

**Partial Summation Formula**
*Lemma* Let $f \in C^{(1)}[0, \infty]$ be a continuously differentiable function. Then

$$\sum_{n \leq x} f(n) = f(1) - (x-[x])f(x) + \int_1^x f(t)\,dt + \int_1^x (t-[t])\,df(t). \tag{79}$$

Proof: Assume $x \in \mathbb{N}$ and use integration by part to link up the finite sum and the integral of the function. The finite sum can be rewritten as

$$\sum_{n \leq x} f(n) = xf(x) - \sum_{n=2}^{x} (n-1)(f(n) - f(n-1)) = xf(x) - \sum_{n=2}^{x} (n-1) \int_{n-1}^{n} df(t)$$
$$= xf(x) - \sum_{n=2}^{x} \int_{n-1}^{n} [t]\,df(t) = xf(x) - \int_1^x [t]\,df(t). \tag{80}$$

On the other hand the finite integral is (by partial integration)

$$\int_1^x f(t)\,dt = xf(x) - f(1) - \int_1^x t\,df(t). \tag{81}$$

**Example 1.** For $x > 0$, $\sum_{n \leq x} \log n = x\log x - x + O(\log x)$.

To see this, take $f(n) = \log x$. Applying the partial summation formula returns

$$\sum_{n \leq x} \log n = \log 1 - (x-[x])\log x + \int_1^x [\log(t)]\,dt + \int_1^x (t-[t])\,d\log t$$

$$= x\log x - x + \int_1^x (t-[t])t^{-1}\,dt = x\log x - x + O(\log x).$$

**Example 2.** Suppose that the Riemann hypothesis is valid, then $\sum_{n \leq x} \frac{\log p}{p} = \log x + 1 + O(x^{-1/2} \log^2 x)$.

To see this, take $f(n) = \log p$ if $n = p$ prime, else $f(n) = 0$, and take $g(n) = 1/n$. Applying the Abel summation formula returns

$$\sum_{n \leq x} \frac{\log p}{p} = S(x)g(x) - \int_x^\infty S(t) dg(t) = \frac{\vartheta(x)}{x} - \int_x^\infty (t + O(t^{1/2} \log^2 t))(-t^{-2}) dt = \log x + 1 + O(x^{-1/2} \log^2 x),$$

where $S(x) = \vartheta(x)$.

**Mellin Transform.** The Mellin transform of a function $f(x)$ and its inverse are given by

$$\breve{f}(s) = \int_0^\infty f(x) x^s \frac{dx}{x} \quad \text{and} \quad f(x) = \frac{1}{i 2\pi} \int_{c-i\infty}^{c+i\infty} \breve{f}(s) x^{-s} ds$$

respectively.

**Fourier Transform.** The Fourier transform of a function $f(x)$ and its inverse are given by

$$\hat{f}(s) = \int_{-\infty}^\infty f(t) e^{i 2\pi s t} dt \quad \text{and} \quad f(t) = \frac{1}{i 2\pi} \int_{-\infty}^\infty \hat{f}(s) e^{-i 2\pi s t} ds$$

respectively. The Guassian density $f(t) = e^{-t^2/2}$ or $f(t) = e^{-\pi t^2}$ is the fixed point of the Fourier transform.

**Poisson Summation Formula.** For a differentiable function $f: \mathbb{R} \to \mathbb{C}$ and its Fourier transform $\hat{f}(s)$, the infinite sums

$$\sum_{-\infty \leq n \leq \infty} f(n) = \sum_{-\infty \leq n \leq \infty} \hat{f}(n)$$

are absolutely convergent. This works just as well in $n$-dimensional space $\mathbb{R}^n$ or $\mathbb{C}^n$.

**Gamma Function.** The gamma function is a generalization of the factorial function on the set of integers to the set of all complex numbers. This is defined by

$$\Gamma(s) = \int_0^\infty x^s e^{-x} \frac{dx}{x},$$

for Re(s) > 0. It is in fact the Mellin transform of the exponential function $f(x) = e^{-x}$. The Weierstrass product expansion

$$\Gamma(s) = s^{-1} e^{-\gamma s} \prod_{n=1}^\infty (1 + s/n)^{-1} e^{s/n}$$

explicitly shows the zeros and poles of this function.

**Recursive Formula**          **1/2-Integer Formula**

$\Gamma(s+1) = s\Gamma(s)$,          $\Gamma(n/2 + 1) = 2^{-(n+1)/2} \pi^{1/2} n!!$, $n \geq 3$.

**Gamma Duplication Formula**
$\Gamma(s)\Gamma(s+1/2) = 2^{1-2s}\pi^{1/2}\Gamma(2s),$

**Gamma Reflection Formula**
$\Gamma(s)\Gamma(1-s) = \pi/\sin(\pi s).$

**Gamma Zeta Integral**
$$\Gamma(s)\zeta(s) = \int_0^\infty \frac{x^{s-1}}{e^x - 1}dx$$

**The Product Rule.** The product of a pair of gamma functions, called the Beta function, is the expression

$$\Gamma(s)\Gamma(t) = \Gamma(s+t)\int_0^\infty x^{s-1}(1-x)^{t-1}dx$$

These formulas are used in the investigation of the analytic continuations, functional equations and properties of its zeros of the zeta function and $L$-functions.

**Asymptotic Expansion Of Log $\Gamma(s)$.** The asymptotic expansion of the log $\Gamma(s)$ is the expression

$$\log\Gamma(s) = (s+1/2)\log s - s + \frac{\log(2\pi)}{2} + \frac{B_2}{2s} + \frac{B_4}{3\cdot 4s^3} + \frac{B_6}{5\cdot 6s^5} + \frac{B_8}{7\cdot 8x^7} + \cdots$$

For a real variable, this is

$$\log x! = (x+1/2)\log x - x + \frac{\log(2\pi)}{2} + \frac{1}{12x} - \frac{1}{360x^3} + \frac{1}{1260x^5} - \frac{1}{1680x^7} + \cdots.$$

The Sterling's Formula is inverse expression

$$N! = e^{-N}N^N\sqrt{2\pi N}\left(1 + \frac{1}{12N} + \frac{1}{288N^2} - \frac{139}{51840N^3} - \right)$$

The inequality $e^{-N+1/(12N+1)}N^N\sqrt{2\pi N} < N! < e^{-N+1/(12N)}N^N\sqrt{2\pi N}$ is also useful in applications.

***p*-adic Sterling's Formula.** Given a pair of integers $x \geq y \geq 0$, the *p*-adic valuation of the binomial coefficient is

$$\text{ord}_p\binom{x+y}{x} = \text{ord}_p((x+y)!) - \text{ord}_p(x!) - \text{ord}_p(y!) = \frac{w_p(x) + w_p(y) - w_p(x+y)}{p-1}$$

where $x = \sum_{0 \leq i \leq k} x_i p^i$, $w_p(x) = \sum_{0 \leq i \leq k} x_i$, and $\text{ord}_p(x!) = x - w_p(x)$.

*Lemma.* Let $p \geq 2$ be a prime. The valuation $v_p(x) = \text{ord}_p(x)$ satisfies the followings.
(i) The maximum prime power $p^v$ divisor of $n!$ is $v_p(n!) = \sum_{i \geq 1}[n/p^i] = [n/p] + [n/p^2] + \cdots$.
(ii) $v_p(n!) < [n/p]\log n/\log p$.         (iii) $v_p(n!/(k!(n-k)!)) < \log n/\log p$.

Proof: For (i) observe that there are $[n/p]$ multiple $a_1 p \leq n$ of $p$ that divides $n!$, there are $[n/p^2]$ multiple $a_2 p^2 \leq n$ of $p^2$ that divides $n!$, …, and so on. And For (ii) notice that
$v_p(n!/(k!(n-k)!)) = \sum_{i=1}^\infty [n/p^i] - [k/p^i] - [(n-k)/p^i] = \sum_{d=1}^\infty (1+(-1)^d)/2 < \log n/\log p$, since $[x+y] - [x] - [y]$ = 0, or 1 for any real numbers $x, y \geq 0$.

**Logarithm Integral.** The first formula

$$li(x) = \int_0^x \frac{dt}{\ln t} = \lim_{\varepsilon \to 0}\left(\int_0^{1-\varepsilon} \frac{dt}{\ln t} + \int_{1+\varepsilon}^x \frac{dt}{\ln t}\right) = c + \int_2^x \frac{dt}{\ln t}, \quad 1 < c < 2,$$

is for general work, and a second formula suitable for numerical calculations is given by

$$li(x) = \gamma + \log\log x + x^{1/2} \sum_{n=1}^{\infty} \frac{(-1)^{n-1}}{n!\, 2^{n-1}} \sum_{k=0}^{(n-1)/2} \frac{1}{2k+1} \log^n x,$$

For a proof, see Ramanujan's Lost Notebook, IV, p. 130.

**Residue Theorem.**

***Residue Theorem*** (Cauchy 1826). Let $f(z)$ be an analytic function on a domain $D \subset \mathbb{C}$, and let $C$ be a closed smooth curve on $D$. Then

$$\int_C f(z)\,dz = i2\pi \sum_{z_i} \mathrm{Re}\,s(f, z_i),$$

where $f(z)$ has its poles at $z_i$, and $\mathrm{Re}\,s(f, z_i) = \lim_{z \to z_i}(z - z_i)f(z)$ is the residue of $f(z)$ at $z_i$.

**Perron Integral.** The Cauchy integral of the function $f(s) = x^s/s$ over the rectangle $a - iT, a + iT, -b - iT, -b + iT$ is computed as four line integrals

$$\frac{1}{i2\pi}\int_R \frac{x^s}{s}\,ds = \frac{1}{i2\pi}\int_{a-iT}^{a+iT} x^s + \frac{1}{i2\pi}\int_{a+iT}^{-b+iT} x^s \frac{ds}{s} + \frac{1}{i2\pi}\int_{-b+iT}^{-b-iT} x^s + \frac{1}{i2\pi}\int_{-b-iT}^{a-iT} x^s = 1.$$

The last equality is the residue $\mathrm{Re}\,s(f, 0) = \lim_{s \to 0}(s - 0)f(s) = 1$ of the function $f(s) = x^s/s$ at the pole at $s = 0$. After simplification, it reduces to the Perron integral

$$\frac{1}{i2\pi}\int_{a-iT}^{a+iT} x^s \frac{ds}{s} = \begin{cases} c\dfrac{x^a}{T\log x}, & 0 < x < 1, \\ 1/2, & x = 1, \\ 1 + c\dfrac{x^a}{T\log x}, & x > 1, \end{cases} \quad (82)$$

where $c > 0$ is a constant.

**Explicit Formula.** An explicit formula for the sum $\sum_{n\geq 1}\Lambda(n)f(n)$ over the primes is a transformation to a sum over the zeros of the zeta function.

***Theorem.*** Let $f$ be a continuously differentiable function on the real line $\mathbb{R}^+$, then

$$\sum_{n\geq 1}\Lambda(n)f(n) = \tilde{f}(0) + \tilde{f}(1) - \sum_{n\geq 1}\tilde{f}(-2n) - \sum_{\rho\neq 0,1}\tilde{f}(\rho), \quad (83)$$

where $\tilde{f}(s)$ is the Mellin transform of $f(x)$, and the index $\rho$ ranges trough the zeros of the zeta function.

Proof: For $a = 1$, and assuming convergence of the sum and integral, one has

$$\sum_{n \geq a} \Lambda(n) f(n) = \sum_{n \geq a} \Lambda(n) \left( \frac{1}{i2\pi} \int_{a-i\infty}^{a+i\infty} \tilde{f}(s) n^{-s} ds \right) = \frac{1}{i2\pi} \int_{a-i\infty}^{a+i\infty} \tilde{f}(s) \sum_{n \geq a} \Lambda(n) n^{-s} ds$$

$$= \frac{1}{i2\pi} \int_{a-i\infty}^{a+i\infty} -\frac{\zeta'(s)}{\zeta(s)} \tilde{f}(s) ds \qquad (84)$$

$$= \tilde{f}(0) + \tilde{f}(1) - \sum_{n \geq 1} \tilde{f}(-2n) - \sum_{\rho \neq 0,1} \tilde{f}(\rho).$$

Another typical application of the Perron formula is in the calculation of the finite sum

$$\sum_{n \leq x} d_k(n) = \frac{1}{i2\pi} \int_{a-i\infty}^{a+i\infty} \zeta^k(s) \frac{x^s}{s} ds,$$

where $d_k(n)$ is the number of $k$-tuples such that $n = a_1 a_2 \cdots a_k$.